\documentclass[leqno, 12pt]{amsart}
\usepackage{amssymb,amsmath,amsthm,graphicx,multirow,setspace,url}
\usepackage{amscd}
\usepackage{tikz-cd}

\usepackage{graphicx}

\usetikzlibrary{arrows}
\usetikzlibrary{quotes}

\title{Local Cohomology of Module of Differentials of integral extensions II}
\author{S. P. Dutta}
\address{Department of Mathematics\\
University of Illinois at Urbana-Champaign\\
1409 West Green Street\\
Urbana, Illinois 61801}

\begin{document}

\newcommand{\hookuparrow}{\mathrel{\rotatebox[origin=c]{90}{$\hookrightarrow$}}}
\newcommand{\hookdownarrow}{\mathrel{\rotatebox[origin=c]{-90}{$\hookrightarrow$}}}

%\numberwithin{equation}{section}
%\numberwithin{table}{section}

%\theoremstyle{custom}
\theoremstyle{plain}
\newtheorem*{theorem*}{Theorem}
\newtheorem{theorem}{Theorem}

\newtheorem{innercustomthm}{Theorem}
\newenvironment{customthm}[1]
  {\renewcommand\theinnercustomthm{#1}\innercustomthm}
  {\endinnercustomthm}

%newtheorem*{theorem}{Theorem}
%
%\newtheoremstyle{named}{}{}{\itshape}{}{\bfseries}{.}{.5em}{\thmnote{#3's }#1}
%\theoremstyle{named}
%\newtheorem*{namedtheorem}{Theorem}

%\newtheorem{theorem}{Theorem}[section]
%\newtheorem*{theorem*}{Theorem}[section]
%\newtheorem{theorem}{Theorem}

%\begin{theorem}  A numbered theorem.    \end{theorem}
%\begin{theorem*} An unnumbered theorem. \end{theorem*}

\newtheorem*{corollary}{Corollary}

\theoremstyle{remark}
\newtheorem*{remark}{Remark}
\newtheorem*{claim}{Claim}

\newtheorem*{namedtheorem}{\theoremname}
\newcommand{\theoremname}{testing}
\newenvironment{named}[1]{\renewcommand{\theoremname}{#1}
  \begin{namedtheorem}}
	{\end{namedtheorem}}
%New function names and their shortcuts

\newcommand{\gr}{\operatorname{grade}}
\newcommand{\Syz}{\operatorname{Syz}}
\newcommand{\syz}[2]{\Syz^{#1}(#2)}
\newcommand{\ring}[2]{{\mathcal{O}_{#1}(#2)}}
\newcommand{\hm}[3]{H_{#1}^{#2}(#3)}
\newcommand{\Tor}{\operatorname{Tor}}
\newcommand{\tor}[3]{\Tor_{#1}^{#2}(#3)}
\newcommand{\Ext}{\operatorname{Ext}}
\newcommand{\ext}[3]{\Ext_{#1}^{#2}(#3)}
\newcommand{\Id}{\operatorname{Id}}
\newcommand{\im}{\operatorname{im}}
\newcommand{\coker}{\operatorname{coker}}
\newcommand{\grade}{\operatorname{grade}}
\newcommand{\Ht}{\operatorname{height}}
\newcommand{\Hom}{\operatorname{Hom}}
\newcommand{\Der}{\operatorname{Der}}
\newcommand{\dm}{\operatorname{dim}}
\renewcommand{\H}{\text{H}}
\newcommand{\Tr}{\text{Tr}}
\newcommand{\D}{\text{D}}
\renewcommand{\L}{\text{L}}
\newcommand{\rank}{\text{rank}}
\renewcommand\rightmark{Local Cohomology of Module of Differentials-II}
\renewcommand\leftmark{S. P. Dutta}

%shortcuts
%\newcommand{\mc}[1]{\mathcal{#1}}
\newcommand{\ul}[1]{\underline{#1}}

\begin{abstract}
In this note ($R, m$) denotes a complete regular local ring and $B$ mostly denotes its absolute integral closure. The four objectives of this paper are the following: i) to determine the highest non-vanishing local cohomology of $\Omega_{B/R}$ in equicharacteristic $0$, ii) to establish a connection between each of $\Omega_{B/R}$ and $\Omega_{B/V}$ and pull-back of $\Omega_{A/V}$ via a short exact sequence together with new observations on corresponding local cohomologies in mixed characteristic where $V$ is the coefficient ring of $R$ and $A$ is its absolute integral closure, iii) to demonstrate that $\Omega_{B/R}$ can be mapped onto a cohomologically Cohen-Macaulay module and iv) to study torsion-free property for $\Omega_{C/V}$ and $\Omega_{C/k}$ along with their respective completions where $C$ is an integral domain and a module finite extension of $R$. In this connection an extension of Suzuki's theorem on normality of complete intersections to the formal set-up in all characteristics is accomplished.
\end{abstract}

\maketitle

\section{Introduction}\label{Intro}
This note is a sequel to our previous work in ([3]). Here $(R, m)$ denotes a complete regular local ring of dimension $n$, $B$ denotes its absolute integral closure unless otherwise stated. In ([3]) it was pointed out that $\H_m^n(\Omega_{B/R}) =0$ and in mixed characteristic $\H_m^{n-1}(\Omega_{B/R}) \neq 0 \iff$ the direct summand property for integral extensions of regular local rings is valid and hence by Andr\'e's work in ([1]) it follows that $\H_m^{n-1}(\Omega_{B/R}) \neq 0$ (theorem 4.1, [3]). However, we were not sure that the above observation on non-vanishing of local cohomology ([6]) is valid in equicharacteristc $0$ and we raised it as a question along with a couple of other questions in ([3]). In this paper we are mainly concerned with the following: i) to provide affirmative answers to some of the questions raised in ([3]) on highest non-vanishing local cohomology of $\Omega_{B/R}$ with respect to $m$ and with respect to any ideal generated by a part of 

\smallskip

\line(1,0){300}

AMS Subject Classification: 13D02, 13D22, 13D45, 13N05

Key words and phrases: module of differentials, integral extension, formal smoothness, local cohomology, absolute integral closure. 

regular system of parameters of length $n-1$ in equicharacteristic $0$, ii) to establish a connection
   between each of $\Omega_{B/R}$ and $\Omega_{B/V}$ and pull back of $\Omega_{A/V}$ via a short exact sequence together with several new observations on corresponding local cohomologies in mixed characteristic where $V$ is the coefficient ring of $R$ and $A$ is the absolute integral closure of $V$, iii) to demonstrate that in mixed characteristic $\Omega_{B/R}$ can be mapped onto  cohomologically Cohen-Macaulay modules $N$ such that for every $t>0$, $(0:p^t)N$ are finitely generated Cohen-Macaulay free $B/p^tB$-modules of rank $n-1$ and dimension $n-1$ and iv) to study torsion free property for standard generators of $\Omega_{B/V}$, $\Omega_{B/k}$ and to study the same property for $\Omega_{C/V}, \Omega_{C/k}$ along with their respective completions where $C$ is an integral domain and a module finite extension of $R$. In this connection an extension of Suzuki's theorem on normality of complete intersection equicharacteristic affine domains to the formal set-up in all characteristics is accomplished. Several new results are also discussed. A brief overview of most of the results of this paper is given below. 

In section 2.1 first we derive the following \textbf{corollary} from theorem 3.5 in ([3]): 

\smallskip
Let $R$ be a normal domain containing a field $k$ of
characteristic $0$ with a non-null derivation $D \in \Der_k (R)$. Let $C$ denote the integral closure of $R$ in a finite extension of $Q(R)$ and $\tilde{C} = \{x \in Q(C)| \Tr(xC) \subset R \}$ ([13]). Let $\tilde{D}: \Omega_{C/k} \rightarrow Q(C)$ denote the $C$-linear map induced by the extension of $D$ to $Q(C)$. Then Im$\tilde{D} \subset \tilde{C}$.

We use the above result to prove the following theorem.

\smallskip
\textbf{Theorem 2.2.} Let $(R, m)$ be a complete regular local ring of dimension $n$ in equicharacteristic $0$ and let $B$ denote its absolute integral closure. Let ${x_2, ..., x_n}$ denote a part of a regular system of parameters of $R$ and $\underline{x}$ denote the corresponding 
ideal. Then $\H{_{\underline{x}}^{n-1}}(\Omega_{B/R}) \neq 0$.

\smallskip
As a \textbf{corollary} (2.3) we derive that $\H_m^{n-1}(\Omega_{B/R}) \neq 0$.

The above theorem settles affirmatively question 2a) in equicharacteristic 0 and the corollary settles affirmatively question 1) in (4.4, [3]).

\smallskip
Using the above corollary we are now able to prove the following result in mixed characteristic, 

\smallskip
\textbf{Proposition 2.5.} The surjective map $\H_m^{n-1}(\Omega_{B/R}) \rightarrow \H{_{\underline x}^{n-1}} (\Omega_{B/R})$ is not an isomorphism (here $p, x_2, ..., x_n$ form a regular system of parameters of $R$).

In section 3 first we point out in lemma 3.1 that if $R_i = V[[X_1, ..., X_i]], \\ i \leq n-1, R=R_{n-1}, R_0 = V, B_i =$ absolute integral closure of $R_i$, $A=B_0$ and $B=B_{n-1}$, then each of $\Omega_{B_i/V}$ and $\Omega_{B_i/R_i}$ is a $B_i$-summand of $\Omega_{B/V}$ and $\Omega_{B/R}$ respectively. Similar statement is also valid in the equicharacteristic case.

For lemma 3.2 and corollaries of theorem 3.3 in 3.4, we assume that [Q(V): Q($\hat{\mathbb{Z}}_{p\mathbb{Z}}$)]$<\infty$

Let $\underline{X}$ denote the ideal generated by $X_1, ..., X_{n-1}$ in $R$.
\\
\textbf{Lemma 3.2} points out the following:

 $\H{_{\underline{X}}^{n-1}}(\Omega_{A/V} \otimes B) \simeq \H_m^{n-1}(\Omega_{A/V} \otimes B) \simeq \H_m^n(B)$.

\smallskip
Thus the validity of the direct summand property for integral extensions ([7]) is equivalent to the non-null property for $\H{_{\underline{X}}^{n-1}}(\Omega_{A/V} \otimes B)$.

\smallskip
In our next theorem we prove the following.

\textbf{Theorem 3.3.} The following sequence is exact: 

\medspace
\ \ \ \ \ \ $0\rightarrow \Omega_{A/V} \otimes B \rightarrow \Omega_{B/V} \rightarrow \Omega_{B/A} \rightarrow 0$

\medspace
In one of the two \textbf{corollaries} of this theorem we point out in 3.4 that

\smallskip
$H_p^1(B)$ imbeds into both $\Omega_{B/V}$ and $\Omega_{B/R}$. 

\smallskip
This implies (Remark 3.4) that both $\Omega_{B/V}$ and $\Omega_{B/R}$ contain isomorphic copies of a Cohen-Macaulay torsion $B$-module. It follows from an exact sequence in proposition 2.5 or in theorem 4.1, ([3]) that both $\Omega_{B/V}$ and $\Omega_{B/k}$ also contain cohomologically Cohen-Macaulay flat $B$-modules. For the later observation no restriction on $Q(V)$ is necessary.

\smallskip
The final theorem of this section demonstrates the following. 

\smallskip
\textbf{Theorem 3.5.} $\Omega_{B/R}$ can be mapped onto a $B$-module $M$ such that i) $M$ is cohomologically Cohen-Macaulay of cohomological dimension $n-1$, ii) for every $t >0, (0:p^t)M$ is finitely generated Cohen-Macaulay free $B/p^tB$-module of rank $n-1$ and dimension $n-1$ and iii) $(0:p^t)M$ is a summand of $(0:p^t) \Omega_{B/R}$.

\smallskip
As a corollary we derive that, for every $t>0$, $\Omega_{B/R}$ contains a finitely generated Cohen-Macaulay free $B/p^tB$-module of rank $n-1$ and dimension $n-1$. Moreover, if [Q(V): Q($\hat{\mathbb{Z}}_{p\mathbb{Z}}$)]$<\infty$, the corresponding rank of this submodule is $n$. 

\smallskip
In section 4 we study torsion-free property for standard generators of each of $\Omega_{B/V}$, $\Omega_{B/k}$ and for each of $\Omega_{C/V}$, $\Omega_{C/k}$ where $C$ is the integral closure of $R$ in a finite extension of $Q(R)$. First we prove the following proposition.

\smallskip
\textbf{Proposition 4.1.} Let $b (\neq 0) \in B$; consider an element $w$ of the form $cdb \in \Omega_{B/V}, c \in B$. If $w$ is a torsion element in $\Omega_{B/V}$, then $b \in A$. 

\ \ \ \ \ \ If $R$ is equicharacteristic i.e. $R=k[[X_1, ..., X_n]]$, then for every $b \in B-\bar{k}$, $db \in \Omega_{B/k}$ is torsion free over $B$ where $\bar{k}$ is the algebraic closure of $k$. Same arguments work in the affine case.

\ \ \ \ \ \ The above proposition induced me to inquire about conditions under which each of $\Omega_{B/V}$ and $\Omega_{B/k}$ would become torsion free as $B$-modules. This led me to Suzuki's theorem ([10]) (statement in (4.3)). Theorem 4.3 extends Suzuki's work on the connection between normality of complete intersection affine domains over a field $k$ and torsion-freeness of the corresponding module of differentials to the formal set-up in all characteristics. To this end we first extend Grothendieck's result: Proposition 22.7.2 in ([5]) to mixed characteristic in the following way: 

\smallskip
\textbf{Theorem 4.2.} Let $V$ be a complete discrete valuation ring in mixed characteristic $p$ with maximal ideal $pV$ such that the residue field $k$ of $V$ is perfect. Let $R = V[[X_1, ..., X_{n-1}]]$ and $C=R/I$ be an integral domain such that $p \neq 0$ in $C$. Let $P$ be a prime ideal of $R$ containing $I$ and let $q=P/I$. The following conditions are equivalent.

i) $C_q$ is formally smooth over $V$ for the q-adic topology. 

ii) There exists derivations $D_i, 1 \leq i \leq r, \in \Der_{_V}(R)$ and elements $\{f_i\}_{1 \leq i \leq r}$ of $I$ such that $\{f_i\}_{1 \leq i \leq r}$ generate $IR_P$ and det$(D_i f_j) \notin P$. 

iii) $C_q$ is an unramified regular local ring.

\smallskip
Next we prove the following extension of Suzuki's theorem.

\smallskip
\textbf{Theorem 4.3.} Let $C=R/I$ where $(R, m)$ is either unramified or equicharacteristic complete regular local ring with coefficient ring $V$(field $k$) such that $V/pV(k)$ is perfect and $I$ is a complete intersection prime ideal in $R$. Let $\widehat \Omega_{C/V}$,  $\widehat \Omega_{C/k}$ denote the $m$-adic completions of $\Omega_{C/V}$ and $\Omega_{C/k}$ respectively. We have the following

1) Suppose that if $p \neq 0 \in C$, every prime ideal $P$ of height $1$ in $C$ containing $p$ is such that $p \notin P^2C_P$. If $C$ is normal then both $\widehat \Omega_{C/V}$ and $\widehat \Omega_{C/k}$ are torsion-free over $C$.

2) If $\widehat{\Omega}_{C/V}$ or $\widehat{\Omega}_{C/k}$ is torsion free then $C$ is normal.

\smallskip
As a \textbf{corollary} we point out that, with $R$ as above, if $C=R[[Y_1, ..., \\ Y_r]]/I$ is a module finite extension of $R$ and is a complete intersection domain such that $Q(C)$ is a separable extension of $Q(R)$, then part 1) of the above theorem implies that $\Omega_{C/V}$ and $\Omega_{C/k}$ are torsion-free over $C$.

\smallskip
Our final \textbf{theorem} (4.5) points out the connection between derivations of $R$ and the torsion sub-module $T$ of $\widehat{\Omega}_{C/V}$, where $C$ is as mentioned in the statement of theorem 4.2. We prove the following: 

for any $\omega \neq 0$ in $\widehat{\Omega}_{C/V}$, $\omega \in T \iff \tilde{D}_i (\omega)=0, \tilde{D}_i$ induced by the extension $D_i$ of the standard derivation $\partial/\partial X_i:R \rightarrow R$ for $1 \leq i \leq n-1$ to $Q(C)$. 

\smallskip
For similar observation in the equicharacteristic case see remark 4) in 4.7.

\smallskip
In \textbf{Corollary 3}, 4.6, we again point out that $\Omega_{B/R}$ has an image $N$ such that $N$ is cohomologically Cohen-Macaulay i.e. $\H_m^i(N)=0$ for $i < n-1$, $\H_m^{n-1}(N) \neq 0$ and $\Omega_{B/V}$ can be mapped on to a $B$-module $W$ such that  $\H_m^i(W) = 0$ for $i \geq 0$. $N$ shares the same properties with $M$ as in theorem 3.5 and its corollary. Moreover, for any $x \neq 0 \in m, (0:x)N$ can be imbedded into a cohomologically Cohen-Macaulay flat $B/xB$-module. 

\smallskip
We provide construction of another such module in remark 5), 4.7. 

\smallskip
We are thankful to B. Bhatt for his note (about a couple months back) through which we became aware of his awesome piece of work in [2]. We used some of his main results to draw our conclusions on Cohen-Macaulayness in several places within this paper.

\smallskip
\textbf{Notations.} Given an extension $A \rightarrow B$ of commutative rings, the corresponding module of differentials is denoted by $\Omega_{B/A}$ (instead of $\Omega^1_{B/A}$); d.v.r. stands for discrete valuation ring; for any integral domain (ring) $C, Q(C)$ denotes the field of fractions (quotient ring) of $C$; for any local ring $(R, m), E$ denotes the injective hull of $R/m$ over $R$, $E(k (P))$ denotes the injective hull of $R/P$ over $R$ and for any $R$-module $M, M^{\vee}$ stands $\Hom_R (M, E)$, $M^*$ stands for $\Hom_R (M, R)$. Unless stated otherwise, in a local ring by local cohomology of a module we mean the same with respect to the maximal ideal.

\section*{Section 2.} 

\smallskip
\textbf{2.1.} For the convenience of the reader we state theorem 3.5, ([3]) and add a sentence or two for indicating its proof.

\smallskip
\textbf{Theorem.} (Theorem 3.5, [3]) Let $R$ be a normal domain of dimension $n$ containing a field $k$ of characteristic 0 with a non-null derivation $D \ \in \Der_k (R)$ and let $B$ be its integral closure in an algebraic extension $F$ of $K$. Then $\Hom_R(\Omega_{B/k}, R) \neq 0$. If $(R, m)$ is a complete local normal domain of dimension $n$ with maximal ideal $m$ and a non-null derivation $D \in \Der_k (R)$, then $\H{_m^n}(\Omega_{B/k})\neq 0$.

\smallskip
The proof of this theorem is obtained by proving the following assertion: Let $K$ denote the field of fractions of $R$. Since characteristic of $k=0$, $\D$ can be extended to a derivation $\D$ (same notation): $B \rightarrow F$ and hence to a derivation: $B_P \rightarrow F$ for any prime ideal $P$ of $R$. Let $\Tr: F \rightarrow K$ denote the trace map. We define an $R$-linear map $\L: \Omega_{B/k} \rightarrow K$ in the following way:

For any $\omega \in \Omega_{B/k}, \L(\omega) = \Tr(\tilde{\D}(\omega))$ where $\tilde{\D}:\Omega_{B/k} \rightarrow F$ is the $B$-linear map induced by $\D$; the same prescription works for $\Omega_{{B_P}/k} \rightarrow F$ for any prime ideal $P$ of $R$. Then im$L \subset R$. 

See ([3]) for details.

Next we derive the following corollary.

\smallskip
\textbf{Corollary.} Let $R$ be a normal domain containing a field $k$ of characteristic $0$ with a non-null derivation $D \in \Der_k (R)$. Let $C$ denote the integral closure of $R$ in a finite extension $Q(C)$ of $Q(R)$. Let $\tilde{C} = \{x \in Q(C)| \Tr (xC) \subset R \}$ ([13]). Then $\tilde{D}(\Omega_{C/k}) \subset \tilde{C}$, where the $C$-linear map $\tilde{D}:\Omega_{C/k} \rightarrow Q(C)$ is induced by the extension of $D$ to $Q(C)$.

\smallskip
\textbf{Proof.} Recall in the proof of the above theorem we have shown that if $L=\Tr_\bullet \tilde{D}:\Omega_{C/k} \rightarrow Q(R)$, then $\im L \subset R$. Let $\lambda \in \im \tilde{D}$, i.e. $\lambda = \tilde{D} (w)$ for some $w \in \Omega_{C/k}$. For any $b \in C$, $b \lambda = b \tilde{D} (w) = \tilde{D} (bw)$. Hence by above the theorem, it follows that $\Tr(b \lambda)= \Tr (\tilde{D} (bw)) \in R$; this implies that $\lambda \in \tilde{C}$. Thus the corollary follows. 

\smallskip
\textbf{2.2.} In the following theorem and its corollary we are going to answer affirmatively a question that was raised in ([3]).

\smallskip
\textbf{Theorem.} Let $R=k[[X_1,...,X_n]]$, characteristic of $k=0$. Let $B$ be the integral closure of $R$ in an algebraic extension of $Q(R)$ such that $B$ contains an element $y$ with the property that $y^2 = X_i$, for some $i$, $1 \leq i \leq n$. Let $\underline{X}$ denote the ideal $(X_1, ..., \hat{X}_i, ... X_n) (\hat{-} \equiv $ omitted). Then $\H_{\underline X}^{n-1} (\Omega_{B/R}) \neq 0$.

\smallskip
\textbf{Proof.} Can assume $i=1$. Let $K=Q(R)$ and $F=Q(B)$. Let $D: B \rightarrow F$ denote the extension of $\partial/ \partial X_1 : R\rightarrow R, \ \partial/ \partial X_1 (X_i)=0$ for $i \neq 1$ and $= 1$ for $i = 1$. Let $\tilde{D}:\Omega_{B/k}\rightarrow F$ denote the $B$-linear map induced by $D$. By hypothesis, 

$y^2=X_1 \Rightarrow$ 2ydy$=dX_1 \Rightarrow \tilde{D}$ (dy)=1/2y and 
$$X_1 dy=0 \in \Omega_{B/R} ........................(1).$$

To prove our assertion it will be enough to show that for $t >0, dy (X_2, \\ ..., X_n)^{t-1} \notin (X_2^t, ..., X_n^t) \Omega_{B/R}$.
If possible, let 
$$dy (X_2, ..., X_n)^{t-1}=\Sigma_{i \geq 2}X_i^t w_i, w_i \in \Omega_{B/R}, 2 \leq i \leq n .......................(2)$$
Consider the short exact sequence (theorem 2.3, lemma 2.4, [3]).

$0 \rightarrow \Omega_{R/k} \otimes B \rightarrow \Omega_{B/k} \xrightarrow{\eta} \Omega_{B/R} \rightarrow 0$

Due to lemma in theorem 3.6, [3] we have a split exact sequence 

$0 \rightarrow \underset{s>0}{\cap} m^s\Omega_{R/k} \rightarrow \Omega_{R/k} \rightarrow {\underset{1}{\overset{n}{\oplus}}} R dX_i \rightarrow 0$

Tensoring with $B$ we obtain from above the following split exact sequence 

$0 \rightarrow \cap{_{s>0}} m^s\Omega_{R/k} \otimes B \rightarrow \Omega_{R/k} \otimes B \rightarrow {\underset{1}{\overset{n}{\oplus}}} B dX_i \rightarrow 0$.

\smallskip
Let $\tilde{w_i} \in \Omega_{B/k}$ be such that $\eta (\tilde{w_i}) = w_i$. The following equality then follows from (2):

\smallskip
$dy (X_2, .., X_n)^{t-1} = \Sigma{_{i \geq 2}} X_i^t \tilde{w_i} + \underset{i \geq 1}{\Sigma} b_i dx_i + \lambda, b_i \in B, \\ \lambda \in \underset{s>0}{\cap} m^s\Omega_{R/k} \otimes B..........................................(3)$

Note that $\tilde{D} (\lambda) = 0$. Applying $\tilde{D} : \Omega_{B/k} \rightarrow F$ to (3) we obtain:

\smallskip 
$$\tilde{D}(dy) (X_2, ..., X_n)^{t-1} = \Sigma X_i^t \tilde{D}(\tilde{w}_i) + b_1$$.

From (1), it follows that 
$$(X_2, ..., X_n)^{t-1} = 2y \Sigma X_i^t \tilde{D}(\tilde{w_i}) + 2y b_1 .........................................(4)$$

Let $S=k[[X_1, ..., X_n, Y]]/(Y^2-X_1) \simeq k [[X_2, ..., X_n, Y]]$.

For any finite extension $C$ of $S$ contained in $B$ such that $b_1 \in C$ and $\tilde{w_i} \in \Omega_{C/k}, 2 \leq i \leq n$, we have $\Tr_{Q(C)/Q(R)}(\mu \tilde{D}(\tilde{w_i})| \mu \in C) \subset R$ (theorem 2.1). By corollary in 2.1 we obtain
$\tilde{D} (\tilde{w_i}) \subset \tilde{C}, \ 2 \leq i \leq n$.

Since $\Tr_{Q(C)/Q(R)} = \Tr_{Q(S)/Q(R)} \bullet \Tr_{Q(C)/Q(S)}$, it follows by the above lemma that $\Tr_{Q(C)/Q(S)} (\tilde{D}(\tilde{w_i})) \subset \tilde{S} = R\{1/2y, y/2y\}, 2 \leq i \leq n$, where $\tilde{S}=\{z \in Q(S)|\Tr(zS) \subset R \}$; $\tilde{S}$ is a free $R$-module generated by $\{1/2y, y/2y\}$ ([9], [13]). 

Applying $\Tr_{Q(C)/Q(S)}$ to (4) we obtain (ignoring the units): 

$(X_2,...,X_n)^{t-1} = \Sigma X_i^t 2y \Tr_{Q(C)/Q(S)} (D(\tilde{w_i})) + 2y \Tr_{Q(C)/Q(S)} b_1$.

It follows from above that 
$$(X_2,...,X_n)^{t-1} = \underset{i\geq2}{\Sigma} X_i^t s_i + 2yc_1, s_i, c_1 \in S$$.
This contradicts the monomial property over $S$ ([6]). Hence \\ $\H_{\underline{X}}^{n-1} (\Omega_{B/R}) \neq 0$.

\smallskip
\textbf{2.3. Corollary.} Let $(R, m)$ and $B$ be as above. Then $\H_m^{n-1} (\Omega_{B/R}) \neq 0$.

\smallskip
\textbf{Proof.} We continue with our notations from the above theorem. Let us recall that localization commutes with module of differentials even in non-noetherian situations. By corollary 2.1 in ([3]) we have $\H_{\underline{X}}^{n-1} (\Omega_{B/R}[1/X_1]) = 0$. This implies, due to the above theorem and the following exact sequence

$\H_m^{n-1} (\Omega_{B/R}) \rightarrow \H_{\underline{X}}^{n-1} (\Omega_{B/R}) \rightarrow \H_{\underline{X}}^{n-1} (\Omega_{B/R}[1/X_1]) \rightarrow 0$

that $\H_m^{n-1} (\Omega_{B/R}) \neq 0$.

\smallskip
\textbf{Remark.} Since $H_m^n(\Omega_{B/R})=0$ (corollary. 2.1, [3]), the above corollary implies that the highest non-vanishing local cohomology of $\Omega_{B/R}$ occurs at the (n-1)-th stage. A proof of this observation can also be obtained without applying corollary 2.1.

\smallskip
\textbf{2.4.} Due to corollary 2 lemma 2.1 in [3] we know if $(R, m)$ is a complete local normal domain of dimension $n$ of equicharacteristic $0$ and $B$ is an integral domain integral over $R$, then $\H_m^n (\Omega_{B/R}) = 0$. In this perspective we have the following proposition.

\smallskip
\textbf{Proposition.}  Let $(R, m)$ be a complete regular local ring of dimension $n$ of equicharacteristic $0$ and $B$ be an integral domain integral over $R$. Suppose $B$ contains a square root of a regular parameter of $R$. Then for any system of parameters $y_1, ..., y_n$ of $R$, $\Omega_{B/R} \otimes B/y\underline{B} = \Omega_{(B/y\underline{B})/(R/ \underline{y})} \neq 0$; here $\underline{y}$ denotes the ideal generated by $y_1, ..., y_n$.

\smallskip
\textbf{Proof.} Let $R=k[[X_1, ..., X_n]]$ and let $\underline{X} = (X_1, ..., X_n)$. It is enough to show that $\Omega_{(B/ \underline{X}B)/k} \neq0$. 
Let $y \in B$ be such that $y^2 = X_1$; then 

$2ydy = dX_1$ i.e. $2X_1 dy = ydX_1$ in $\Omega_{B/k}$. 

If possible let $\im dy = 0$ in $\Omega_{(B/ \underline{X})/k}$. 

Consider the exact sequence (theorem 2.3, lemma 2.4, [3]): 
$$0 \rightarrow \Omega_{R/k} \otimes B \rightarrow \Omega_{B/k} \rightarrow \Omega_{B/R} \rightarrow 0$$
Tensoring the sequence with $B/ \underline{X}B$ we obtain the following exact sequence: 
$$\rightarrow \Omega_{R/k} \otimes B/ \underline{X}B \rightarrow \Omega_{B/k} \otimes B/\underline{X}B \rightarrow \Omega_{(B/ \underline{X}B)/k} \rightarrow 0$$ 
im$dy=0$ in $\Omega_{(B/ \underline{X}B)/k}\Rightarrow$ 

$dy = \Sigma X_i w_i + \Sigma b_i dX_i + \lambda$ in $\Omega_{B/k}$, $b_i \in B$, $w_i \in \Omega_{B/k}$ for $1 \leq i \leq n$ and $\lambda \in \cap m^t \Omega_{R/k} \otimes B$.

Let $\tilde{D}: \Omega_{B/k} \rightarrow Q(B)$ denote the $B$-linear map induced by the extension of $\partial/ \partial X_1: R \rightarrow R$.
Applying $\tilde{D}$ to the above equation, we obtain 

$$1/2y = \Sigma X_i D(w_i) + b_1.$$ 

This implies 
$$1 - 2yb_1 = \Sigma X_i 2y \tilde{D} (w_i).$$
Since $R$ is complete, $B$ is local and $1-2y b_1$ is unit, say $u$ in $B$. Thus 
$$1=\Sigma X_i 2i u^{-1} \tilde{D}(w_i).$$
Applying the trace map we obtain via theorem 2.1.
 
$1 = \Sigma X_i r_i, r_i \in R$ --- a contradiction.

\smallskip
\textbf{Remark.} This result is valid even if $R$ is not complete provided $B$ is local.

\smallskip
\textbf{2.5.} Here we will use theorem 2.3 to prove an observation in the mixed characteristic.
 
\smallskip
\textbf{Proposition 2.5.} Let $R = V[[X_1, ..., X_{n-1}]]$, $m = (p, X_1, ..., X_{n-1})$ and $(V, pV)$ is a complete discrete valuation ring. let $B$ denote the absolute integral closure of $R$ and let $\underline{X}$ denote the ideal generated by $X_1, ..., X_{n-1}$ in $R$. Then the map $\alpha : \H_{\underline{X}}^{n-2} (\Omega_{B/R}[1/p]) \rightarrow \H_m^{n-1} (\Omega_{B/R})$ is non-null.

\smallskip
\textbf{Proof.} Let us recall that the above map is part of a long exact sequence: 
$$\H_{\underline{X}}^{n-2} (\Omega_{B/R}) \rightarrow \H_{\underline{X}}^{n-2} (\Omega_{B/R}[1/p]) \xrightarrow{^\alpha} \H_m^{n-1} (\Omega_{B/R}) \rightarrow \H_{\underline{X}}^{n-1} (\Omega_{B/R}) .......(5)$$

Let $y\in B$ be such that $y^p=X_1 \Rightarrow p y^{p-1} dy=dX_1 \Rightarrow pX_1 dy=ydX_1 \in \Omega_{B/V} \Rightarrow pX_1 dy = 0 \in \Omega_{B/R}$.

Consider the element $\lambda =(dy(X_2, ..., X_{n-1})^{t-1}, 0 ... 0) \in \oplus{_1 ^{n-1}} \Omega_{B/R}[1/p]$.

Note that $X_1^t dy(X_2, ..., X_{n-1})^{t-1} = X_1 dy (X_1, ..., X_{n-1})^{t-1} = 0 \in \Omega_{B/R}[1/p]$.
Thus $\lambda$ represents a non-null element in $\H^{n-1}(\underline{X}^t, \Omega_{B/R}[1/p])$ (proof of theorem 2.2). Since $\Omega_{B/R}(\Omega_{B/V})=p \Omega_{B/R}(\Omega_{B/V})$ (theorem 4.1, [3]), we have an exact sequence 
$$0 \rightarrow \H_p^0 (\Omega_{B/R}) \rightarrow \Omega_{B/R} \rightarrow \Omega_{B/R} [1/p] \rightarrow 0$$.
It follows that $\H_{\underline{X}}^i (\H_p^0 (\Omega_{B/R})) \simeq \H_m^i (\Omega_{B/R})$ for $i \geq 0$. Then $\alpha ($class of $\lambda)=$ class of $X_1 dy (X_1, ..., X_{n-1})^{t-1} \in \H_{\underline{X}}^{n-1} (\H_p^o (\Omega_{B/R}))$. 
So $\alpha($ class of $\lambda)=0 \Rightarrow$:
$$X_1 dy (X_1, ..., X_{n-1})^{t-1} = \Sigma{_{i \geq 2}} X_i^t w_i, \ \ w_i \in (0:p^t) \Omega_{B/R}, 2\leq i\leq{n-1} ....(6)$$
We have an exact sequence (theorem 2.3, lemma 2.4, [3]): 
$$0 \rightarrow \Omega_{R/V} \otimes B \rightarrow \Omega_{B/V} \rightarrow \Omega_{B/R} \rightarrow 0$$
Applying $\Hom_B (B/p^t B, -)$ we obtain the following short exact sequence:
$$0 \rightarrow (0:p^t) \Omega_{B/V} \rightarrow (0:p^t) \Omega_{B/R} \xrightarrow{^\beta} \Omega_{R/V} \otimes B/p^tB \rightarrow 0 .....(7)$$
Applying $\beta$ to (6) we obtain the following equality from equ (7):

$p^{t-1}ydX_1(X_1, ..., X_{n-1})^{t-1}=\underset{i\geq 1}{\Sigma} X_i^t \mu_i \in \Omega_{R/V} \otimes B/p^tB, \ \mu_i \in \Omega_{R/V} \otimes B, 1 \leq i \leq n$.

Hence, 
$$p^{t-1}ydX_1(X_1, ..., X_{n-1})^{t-1}=\underset{i\geq 1}{\Sigma} X_i^t \mu_i + p^t \nu, \nu \in \Omega_{R/V} \otimes B ........(8)$$

Let $D_1$ be the extension of $\partial/ \partial X_1: R \rightarrow R$ to $B \rightarrow Q(B)$ and let $\tilde{D_1}: \Omega_{B/k} \rightarrow Q(B)$ denote the corresponding B-linear map.

Applying $\tilde{D_1}$ to (8) we obtain

\smallskip
$p^{t-1} y(X_1, ..., X_{n-1})^{t-1}=\Sigma X_i^t b_i + p^t c$ where $b_i, c \in B, 1 \leq i \leq n$.

\smallskip
Since $y^p=X_1$, the above equality contradicts the monomial property in $B$. Hence im$\alpha \neq 0$.

\smallskip
\textbf{Remark.} Since  $\H_{\underline{X}}^{n-1}(\Omega_{B/R}[1/p]) = 0$ (corollary 2.1, [3]), the above proposition shows that the surjective map $\H_m^{n-1}(\Omega_{B/R}) \rightarrow \H_{\underline{X}}^{n-1}(\Omega_{B/R})$ is not an isomorphism. Let us mention that we are still not able to prove that $\H_{\underline{X}}^{n-1}(\Omega_{B/R}) \neq 0$ in mixed characteristic (raised as a question in [3]).

\smallskip
\section*{Section 3.}

\smallskip
\textbf{3.1. Lemma.} Let $R_i = V[[X_1, ..., X_i]]$, $1\leq i \leq n-1$, $R_0 =V$ and $R_{n-1} = R$. Let $B, B_i$ denote the absolute integral closure of $R, R_i$ respectively for $0 \leq i \leq n-1$, $B_0 = A$. Then each of $\Omega_{B_i/V}$ and $\Omega_{B_i/R_i}$ is $B_i$-summand of $\Omega_{B/V}$ and $\Omega_{B/R}$ respectively for $0 \leq i <n-1$. 

Similar result is also valid in equicharacteristic zero. 

\smallskip
\textbf{Proof.} Let $Q_i$ be a prime ideal of $B$ lying over $(X_{i+1}, ..., X_{n-1})$ in $R$. Since $B$ is absolutely integrally closed, we have $Q_i = Q_i^2$, $B/Q_i$ is the absolute integral closure of $R_i$ and hence $B_i \simeq B/Q_i$. Thus $B_i \xrightarrow{J_i} B \xrightarrow{\eta_i} B/Q_i$ is an $R_i$-algebra isomorphism.
We have the following exact sequences : 
$$Q_i/Q_i^2 \rightarrow \Omega_{B/V} \otimes B/Q_i \rightarrow \Omega_{(B/Q_i)/V} \rightarrow 0$$ and
$$Q_i/Q_i^2 \rightarrow \Omega_{B/R} \otimes B/Q_i \rightarrow \Omega_{(B/Q_i)/R_i} \rightarrow 0.$$ 

it follows from above that 
$$\Omega_{B/V} \otimes B/Q_i \simeq \Omega_{B_i/V}$$ and
$$\Omega_{B/R} \otimes B/Q_i \simeq \Omega_{B_i/R_i}.$$

Thus the composite maps:
$$\Omega_{B_i/V} \xrightarrow{\tilde{J_i}} \Omega_{B/V} \xrightarrow{\tilde{\eta_i}} \Omega_{(B/Q_i)/V}$$ and
$$\Omega_{B_i/R} \xrightarrow{\tilde{J_i}} \Omega_{B/R} \xrightarrow{\tilde{\eta_i}} \Omega_{(B/Q_i)/R_i}$$

are $B_i$-module isomorphism. Hence the result.

\smallskip
\textbf{3.2.}
For the rest of section 3 we assume $R = V[[X_1, ..., X_{n-1}]]$, $V$-a complete discrete valuation ring in mixed characteristic $p$ such that $V/pV$ is perfect; let $A,B$ denote the absolute integral closures of $V, R$ respectively. For lemma 3.2 and corollary 2 in 3.4 we assume that $[Q(V):Q_p] < \infty, Q_p = Q(\hat{Z}_{pZ})$. 

Since $\Omega_{Q(A)/Q(V)} =0, \Omega_{A/V}$ is a $p$-torsion module. In a more general setup Fontaine [4] showed that $\Omega_{A/V}$ is generated by $\{d \xi_t \}_{t\geq1}$, where $\xi_t$ is a $p^t$-th root of unity such that $\xi_{t+1}^p=\xi_t$ for $t>0$. We have the following lemma.
 
\smallskip
\textbf{Lemma.} $\Omega_{A/V} \otimes B \simeq H_p^1 (B)$ and $H_m^{n-1} (\Omega_{A/V} \otimes B) \simeq H_m^n (B)$.

\smallskip
\textbf{Proof.} Since $\Omega_{A/V}$ is a $p$-torsion module, $\Omega_{A/V} \simeq \underrightarrow{\lim} (0:p^t) \Omega_{A/V}$. In ([12]) it was derived from Fontaine's work that $(0:p^t) \Omega_{A/V} \simeq A/p^t A$.
We have the following commutative diagram:

\[
\begin{CD}
	@. @. (0:p^t) \Omega_{A/V} @>>> A/p^t A \\
	@. @. @VVV  @VVV \\
	@. @. (0:p^{t+1}) \Omega_{A/V} @>>> A/p^{t+1} A \\
%	@.  @VVV  @VVV \\
\end{CD}
\]

Where the horizontal arrows are isomorphisms, the right vertical arrow is multiplication by $p$ and the left vertical arrow is the inclusion map.

It follows from above that $\Omega_{A/V} \simeq \underrightarrow{\lim} A/p^t A = \H_p^1 (A)$.

Hence $\Omega_{A/V} \otimes B \simeq \H_p^1 (B)$. 

Moreover $\H_m^{n-1} (\Omega_{A/V} \otimes B) \simeq \H_m^{n-1} (\H_p^1 (B)) \simeq \H_m^n(B)$

\smallskip
\textbf{Corollary.} The validity of the monomial property for local rings is equivalent to the condition that $\H_m^{n-1} (\Omega_{A/V} \otimes B) \neq 0$, i.e. $\H_{\underline{X}}^{n-1} (\Omega_{A/V} \otimes B) \neq 0$

\smallskip
\textbf{3.3.} Our next observation is the following.

\smallskip	
\textbf{Theorem.} Let $AR=A \underset{V}{\otimes}R$. With notations as above the following sequence is exact:

\smallskip
$0 \rightarrow \Omega_{A/V} \underset{A}{\otimes}B \overset{\alpha}{\rightarrow} \Omega_{B/R} \rightarrow \Omega_{B/AR} \rightarrow 0$

\smallskip
\textbf{Proof.} The following sequence is exact
$$\Omega_{A/V} \underset{A}{\otimes} B \overset{\alpha}{\rightarrow} \Omega_{B/R} \rightarrow \Omega_{B/AR} \rightarrow 0.$$
We need to show that $\alpha$ is injective.

We have $A=\underrightarrow{\lim}A_t$ where, for every $t$, $A_t$ is a complete d.v.r. and a module finite extension of $V$;

$A_t=V[Y]/(f_t(Y)), f_t(Y)$ is a monic irreducible polynomial in $V[Y]$ ([9]).

$\Omega_{A/V}=\underrightarrow{\lim} \Omega_{A_t/V}$. It is enough to show that

$\Omega_{A_t/V} \otimes B \rightarrow \Omega_{B/R}$ is injective for every $t$ .............(9).

$B=\underrightarrow{\lim}C$, $C$ is a module finite extension of $R$.
 
$\Omega_{B/R}=\underrightarrow{\lim} \Omega_{C/R}$.
$A_tR=A_t \underset{V}{\otimes}R \simeq A_t[[X_1, ..., X_{n-1}]] \simeq R[Y]/(f_t(Y))$; $A_tR$ is a complete regular local ring.

$\Omega_{A_t/V} \underset{A_t}{\otimes} B \simeq \Omega_{A_t/V} \underset{A_t}{\otimes} A_tR \underset{A_tR}{\otimes} B$;

$\Omega_{A_t/V} \underset{A_t}{\otimes} A_tR \simeq \Omega_{A_t/V} \underset{V}{\otimes}R\simeq \Omega_{A_tR/R}$.

Write $f=f_t(Y)$, $f'=f'_t(Y)$ and $S=A_tR=R[Y]/(f)$. We need to show $\Omega_{S/R} \underset{S}{\otimes}B \rightarrow \Omega_{B/R}$ is injective.

It is enough to show that $\Omega_{S/R} \underset{S}{\otimes}C\rightarrow \Omega_{C/R}$ is injective where $C$ is any normal domain module finite extension of $S$ (hence a module finite extension of $R$). 

$\Omega_{S/R}=S/f'S$; hence $\Omega_{S/R} \underset{S}{\otimes}C=C/f'C$.

We have an exact sequence ([11]): 
$$\Gamma_{C/S} \overset{\beta}{\rightarrow} \Omega_{S/R} \underset{S}{\otimes}C \overset{v}{\rightarrow} \Omega_{C/R} \rightarrow \Omega_{C/S} \rightarrow 0$$
we need to show $v$ is injective i.e. to show that im$\beta=0$.

Since $S$ and $C$ are complete local normal domains and $C$ is a module finite extension of $S$, for every prime ideal $P$ of height$1$ in $S$, $C_P$ is a Dedekind domain module finite extension of $S_P$. Hence, by lemma 2.2 in [3] $(\Gamma_{C/S})_P=\Gamma_{C_P/S_P}=0$. This implies that im$\beta$ has at least codimension $1$ in $\Omega_{S/R} \otimes C=C/f'C$. Since $C$ is a complete local normal domain, $C/f'C$ can't have any sub-module of codimension $\geq1$ in $C/f'C$. Thus im$\beta=0$ and the assertion follows.

\smallskip
\textbf{3.4. Corollaries.}

\smallskip
\textbf{1.} There exists a commutative diagram of short exact sequences 

\[
\begin{CD}
	@. @. 0 @. 0 \\
	@. @. @VVV  @VVV \\
	@. @. \Omega_{A/V} \otimes B @= \Omega_{A/V} \otimes B @. \\
	@. @. @VVV  @VVV \\
	0 @>>>\Omega_{R/V} \otimes B @>>> \Omega_{B/V} @>>> \Omega_{B/R} @>>> 0 \\
	%  @VVV  @. @VVV \\
	@.   \parallel @. @VVV @VVV @. \\
	0 @>>>\Omega_{R/V} \otimes B @>>> \Omega_{B/A} @>>> \Omega_{B/AR} @>>> 0 \\
	@. @. @VVV  @VVV \\
	@. @. 0 @. 0
\end{CD}
\]

\smallskip
\textbf{Proof.} The exactness of the middle row follows from theorem 2.3 and lemma 2.4 in [3]. Since $\Omega_{A/V} \otimes B$ is a $p$-torsion module and $\Omega_{R/V} \otimes B$ is $B$-flat and hence torsion free over $B$, the assertion follows. For local cohomological implications see Remarks 2 in (4.7).

\medskip
\textbf{2.} There exist exact sequences: 
$$0 \rightarrow \H_p^1 (B) \rightarrow \H_p^0 (\Omega_{B/V}) \rightarrow \H_p^0 (\Omega_{B/A}) \rightarrow 0$$ and
$$0 \rightarrow \H_p^1 (B) \rightarrow \H_p^0 (\Omega_{B/R}) \rightarrow \H_p^0 (\Omega_{B/AR}) \rightarrow 0$$.

\smallskip
\textbf{Proof.} By the above corollary we have an exact sequence:
$$0 \rightarrow \Omega_{A/V} \otimes_{A}B\rightarrow \Omega_{B/V} \rightarrow \Omega_{B/A} \rightarrow 0.$$
This implies, due to lemma (3.2), and the fact that $\Omega_{A/V}=p \Omega_{A/V}$, the exactness of the first exact sequence.

%Applying $\Hom_B(B/p^t B, -)$ to the above exact sequence, we obtain for every $t>0$ the following exact sequence:
%$$0 \rightarrow(0:p^t) (\Omega_{A/V}\otimes_AB) \rightarrow(0:p^t) \Omega_{B/V}\rightarrow (0:p^t)\Omega_{B/A} \rightarrow 0.$$
%$\Omega_{A/V}=\underrightarrow{\lim}\Omega_{{A_t}/V}=\underrightarrow{\lim}A_t/(f'(y_t))$, where $A_t$ is a DVR $\subset A, A_t=V[y_t]$, where $f(X)$ denotes the minimal polynomial for $\xi_t$ over $V$, $f'(y_t) \neq 0$ in $A_t$ ([8]). Since $A_t \hookrightarrow B$ is injective and $B$ is a domain, $\Tor{_1^{A_t}}(A_t/(f' (y_t)), B)=0$. Hence $\Tor_1^A(\Omega_{A/V}, B)=0$. Consider the short exact sequence:
%$$0\rightarrow(0:p^t)\Omega_{A/V}\rightarrow\Omega_{A/V}\stackrel{p^t}{\rightarrow}\Omega_{A/V}\rightarrow 0.$$
%Tensoring the above sequence with $B$ we obtain $(0:p^t)(\Omega_{A/V}\otimes_AB)\simeq (0:p^t) \Omega_{A/V} \otimes_AB \simeq A/p^tA \otimes_AB$, $=B/p^tB$. Thus we have the following exact sequence:
%$$0 \rightarrow B/p^t B \rightarrow (0:p^t) \Omega_{B/V} \rightarrow (0:p^t)\Omega_{B/A}\rightarrow 0.$$
%Taking direct limit, we obtain the desired exact sequence. 

The second exact sequence is also proved by similar arguments.

\smallskip
\textbf{Remark.} Corollary 2 above implies, due to corollary 5.10 [2], that both $\Omega_{B/V}$ and $\Omega_{B/R}$ contains a Cohen-Macaulay torsion $B$-module. From middle row of the diagram in corollary 1 above it follows, due to theorem 5.1 [2] and flatness of $\Omega_{R/V}$ over $R$, that $\Omega_{B/V}$ also contains a cohomologically Cohen-Macaulay flat $B$-module. For the second assertion no restriction on $Q(V)$ is necessary. 

\smallskip
\textbf{3.5.} In our final theorem of this section we prove the following.

\smallskip
\textbf{Theorem.} With $V, R, B$ as in 3.2, $\Omega_{B/R}$ can be mapped onto a $B$-module $M$ such that i) $M$ is cohomologically Cohen-Macaulay of cohomological dimension $n-1$, ii) for every $t >0, (0:p^t)M$ is finitely generated Cohen-Macaulay free $B/p^tB$-module of rank $n-1$ and dimension $n-1$ and iii) $(0:p^t)M$ is a summand of $(0:p^t) \Omega_{B/R}$.

\smallskip
\textbf{Proof.} We have a short exact sequence:
$$0 \rightarrow \H_p^0 (\Omega_{B/V}) \rightarrow \Omega_{B/V} \rightarrow \Omega_{B/V} [1/p] \rightarrow 0.$$
(the exactness follows from the fact that $\Omega_{B/V} \otimes B/pB=0$(theorem 4.1, [3]))

Let $M$ be the cokernel of $(\H_p^0 (\Omega_{B/V}) \rightarrow \Omega_{B/R})$

\smallskip
Consider the following commutative diagram of short exact sequences:

\[
\begin{CD}
	@. @. 0 @. 0 \\
	@. @. @VVV  @VVV \\
	@. @. \H_p^0 (\Omega_{B/V}) @= \H_p^0 (\Omega_{B/V}) @. \\
	@. @. @VVV  @VVV \\
	0 @>>>\Omega_{R/V} \otimes B @>>> \Omega_{B/V} @>>> \Omega_{B/R} @>>> 0 \\
	%  @VVV  @. @VVV \\
	@.   \parallel @. @VVV @VVV @. \\
	0 @>>>\Omega_{R/V} \otimes B @>>> \Omega_{B/V}[1/p] @>>> M @>>> 0 \\
	@. @. @VVV  @VVV \\
	@. @. 0 @. 0
\end{CD}
\]

\smallskip
The exactness of the left column follows from the fact that $\Omega_{R/V}$ being a flat $R$-module ([3]), $\Omega_{R/V} \otimes B$ is a flat $B$-module and hence torsion-free over $B$.

\smallskip
Since $\H_m^i(\Omega_{B/V}[1/p])=0$ for $i\geq 0$ and $\Omega_{R/V}$ is $R$-flat (theorem 2.3, [3]), it follows from the last row of the above diagram that $\H_m^i(M) \simeq \Omega_{R/V} \otimes \H_m^{i+1}(B) \simeq {\underset{1} {\overset{n-1}{\oplus}}} \H_m^{i+1}(B)$ (lemma, theorem 3.6, [3]) for $i \geq 0 ............................................................................................(*)$; 

in particular $\H_m^{n-1}(M) \simeq {\underset{1} {\overset{n-1}{\oplus}}} \H_m^n (B)$.

\smallskip
By theorem 5.1, [2] it follows that $\H_m^i(B)=0$ for $i<n$; hence $\H_m^i (M)=0$ for $i<n-1$ and $\H_m^{n-1}(M) \neq 0$ ([1], [7]). Thus $M$ is cohomologically Cohen-Macaulay of cohomological dimension $n-1$. From the last row of the above commutative diagram we obtain $(0:p^t)M \simeq \Omega_{R/V} \otimes B/p^tB \simeq {\underset{1} {\overset{n-1}{\oplus}}} B/p^tB$ (lemma within theorem 3.6, [3]). Since, by corollary 5.10 [2], $B/p^tB$ is Cohen-Macaulay, our assertion follows.

\smallskip
\textbf{Corollary.} For every $t>0$, $\Omega_{B/R}$ contains a finitely generated Cohen-Macaulay free $B/p^tB$-module of rank $n-1$ and dimension $n-1$. If $[Q(V):Q(\hat{Z}_{pZ})]<\infty$, then the corresponding rank is $n$.

Proof follows from: a) the above commutative diagram by applying $\Hom(B/p^tB, -)$ to the last row and the last column of the above diagram, b) corollary 2, 3.4 and c) corollary 5.10, [2].

\smallskip
\textbf{Remark.} Since $\H_m^{n-1}(\Omega_{B/R}) \rightarrow \H_m^{n-1}(M) \rightarrow 0$ is exact, the observation in (*) also leads to a new proof of the fact (theorem 4.1, [3]) that $\H_m^{n-1} (\Omega_{B/R}) \neq 0$ if and only if the direct summand property is valid over $R$.

\smallskip
\section*{Section 4.}
\textbf{4.1.} Our first proposition deals with torsion-free property of standard generators of $\Omega_{B/V}$.

\smallskip
\textbf{Proposition.} Let $R = V[[X_1, ... , X_{n-1}]]$ or $k[[X_1, ... , X_n]]$ or $k[X_1, ... , \\ X_n]$ where characteristic of $k$ is $0$. Let $B$ denote an integral extension of $R$ in an algebraic extension of $Q(R)$. If for some $c \in B, b \in B-V, cdb \in \Omega_{B/V}$ is a torsion element, then $b \in A \cap B$ where $A$ is the absolute integral closure of $V$. 

\smallskip
\ \ \ \ \ \ \ \ If $R$, as above, contains a field $k$ of characteristic $0$ then for $b \in B-\bar{k}, db$ is torsion-free, where $\bar{k}$ is the algebraic closure of $k$ in $Q(B)$.

\smallskip
\textbf{Proof.} Let $R = V[[X_1, ..., X_{n-1}]]$. Let $b \in B$ be such that $db$ is a non-null torsion element in $\Omega_{B/V}$. We have an exact sequence 
$$0 \rightarrow \Omega_{R/V} \otimes B \overset{\beta} {\rightarrow} \Omega_{B/V} \rightarrow \Omega_{B/R} \rightarrow 0$$
Let $c' \in B-(0)$ be such that $c'cdb = 0$. Let $\tilde{c}=c'c; \ \tilde{c}db=0$. $\tilde{c}$ satisfies an integral equation over $R$: $\tilde{c}^s + r'_1 \tilde{c}^{s-1} + ... + r'_s = 0$, $r'_i \in R, 1 \leq i \leq s$ and $r'_s \neq 0$. This implies that $r'_s db = 0$. We write $r$ for $r'_s$ i.e. $rdb =0, r\neq 0 \in R$. 

\smallskip
Since $Q(B)$ is separable algebraic over $Q(R)$, if $f(X)$ denotes the minimal polynomial for $b$ over $R$ then $f'(b) \neq 0$. Hence
$$f'(b)rdb=0..........................................(10)$$
Let $f(X)=X^t+r_1 X^{t-1}+ ... +r_t$. Then 
$$f(b)=b^t+r_1b^{t-1}+ \dots + r_t =0....................(11)$$
Hence $f'(b)db+ \Sigma d(r_i)b^{t-i}=0$ in $\Omega_{B/V}$. This implies, due to (10)
$$r\Sigma d(r_i) b^{t-i} = 0..........................(12)$$
Since $\Omega_{R/V}$ is flat over $R$ and $B$ is an integral domain, $\Omega_{R/V} \otimes_RB$ is torsion-free $B$-module. From the above exact sequence we obtain, by (12),
 
$\Sigma d(r_i) \otimes b^{t-i} = 0$ in $\Omega_{R/V} \otimes_RB$, i.e. $\Sigma d(r_i)b^{t-i}=0 \in \Omega_{B/V}$.

$R[b]$ is a free $R$-module of rank $t$ (due to (10)) with basis $1, b \dots b^{t-1}$. We have an isomorphism
$$\eta: \oplus_1^tR \rightarrow R[b], \ \eta (e_i) = b^{t-i}........................(13)$$
Tensoring (13) with $\Omega_{R/V}$ we get an isomorphism, 
$$\tilde{\eta}: \oplus_1^t \Omega_{R/V} \rightarrow \Omega_{R/V} \otimes_RR[b], \ \tilde{\eta}(d(s_1), \dots, d(s_t))=\Sigma d(s_i) \otimes b^{t-i}.$$ 
Since $\Omega_{R/V}$ is $R$-flat, $\alpha$: $\Omega_{R/V} \otimes_RR[b] \rightarrow \Omega_{R/V} \otimes_RB$ is injective. 

Composing this with the injection $\beta: \Omega_{R/V} \otimes_RB \rightarrow \Omega_{B/V}$ we have an injective map $\beta \cdot \alpha: \Omega_{R/V} \otimes_RR[b] \rightarrow \Omega_{B/V}$. 

Since $\beta \cdot \alpha (\Sigma d(r_i) \otimes b^{t-i} = \Sigma d(r_i) b^{t-i}) = 0$, and $\tilde{\eta}$ is an isomorphism we have $d(r_i) = 0, 1 \leq i \leq t$. This implies that $r_i \in V, 1\leq i \leq t$. 

Hence from (11) it follows that $b$ is integral over $V$ and thus $b \in A$.

If $R$ is of equicharacteristic $0$, the above arguments show that for every $b \in B-\bar{k}, db$ is torsion-free over $B$.

\smallskip
\textbf{4.2.} As mentioned in the introduction proposition 4.1 led us to look for sufficient conditions for torsion-free property of the module of differentials as a whole. And in this respect we were guided by Suzuki's work ([10]). Our main goal in 4.2 and 4.3 is to prove a formal version in all characteristics of Suzuki's theorem on the connection between normality of complete intersection affine domains over a field $k$ and torsion-freeness of the corresponding module of differentials. For this purpose we first extend Grothendieck's proposition 22.7.2 in ([5]) to mixed characteristic. Grothendieck's arguments work in this case with several minor modifications.

\textbf{Theorem.} Let $V$ be a complete discrete valuation ring in mixed characteristic $p$ with maximal ideal $pV$ such that the residue field $k$ of $V$ is perfect. Let $R = V[[X_1, ..., X_{n-1}]]$ and $C=R/I$ be an integral domain such that $p \neq 0$ in $C$. Let $P$ be a prime ideal of $R$ containing $I$ and let $q=P/I$. The following conditions are equivalent.

\smallskip
i) $C_q$ is formally smooth over $V$ for the q-adic topology. 

\smallskip
ii) There exists derivations $D_i, 1 \leq i \leq r, \in \Der_{_V}(R)$ and elements $\{f_i\}_{1 \leq i \leq r}$ of $I$ such that $\{f_i\}_{1 \leq i \leq r}$ generate $IR_P$ and det$(D_i f_j) \notin P$. 

\smallskip
iii) $C_q$ is an unramified regular local ring.

\smallskip
\textbf{Proof.} i) $\iff$ iii). Since $V$ is a d.v.r. and $C_q$ is a domain, $V \rightarrow C_q$ is formally smooth if and only if $V/pV \rightarrow C_q/pC_q$ is formally smooth. Since $V/pV$ is perfect this is equivalent to $C_q/pC_q$ being a regular local ring, i.e. $C_q$ is an unramified regular local ring.

\smallskip
The proof of i) $\iff$ ii) becomes evident by following Grothendieck's arguments in his proof of 22.7.2 in ([5]) with relevant minor modifications needed for the mixed characteristic case. Grothendieck's assumption for ensuring separability is assured here by assuming that the residue field $k$ of $V$ is perfect.

\smallskip
First we recall two basic facts that will be needed in our proof.

\smallskip
a) If $p$ is in $J$ for any prime ideal $J$ of $R$ then $p \notin J^2 R_J$ ([8]).

\smallskip
b) $\widehat{\Omega}_{R/V} \simeq {\overset{n-1}{\underset{1}{\oplus}}} R dX_i, 1 \leq i \leq n-1$, where $\widehat{\Omega}_{R/V}$ is the m-adic completion of $\Omega_{R/V}, m =$ the maximal ideal of $R$ ([3], [5]).

Let $L = k(p)=R_P/P R_P = C_q/q C_q$. By (22.6.2, (ii), [5]) the condition i) is equivalent to the condition 

$d_0:I/I^2 \otimes L \rightarrow \Omega_{R/V} \otimes L$ is injective, $d_0$ is induced by the standard differential map $"d"$. Let us consider the composite map 

$\hat{d}:I/I^2 \otimes L \rightarrow \Omega_{R/V} \otimes L \rightarrow \widehat{\Omega}_{R/V} \otimes L ..............(14)$

and show that condition ii) is equivalent to the assertion that $\hat{d}$ is injective. Note that $I/I^2 \underset{C}{\otimes} L = IR_P/I^2R_P \underset{RP}{\otimes} L$. Condition ii) implies that, if $F_i=\hat{d} (\bar{f_i} \otimes 1)$, where $\bar{f_i}$ is the image of $f_i$ in $IR_P/I^2R_P$, the matrix $(<F_i, D_j \otimes 1>)$ is invertible. This implies $\{F_i\}$ are linearly independent, and hence, so are $\{\bar{f_i} \otimes 1\}$ which generate $ IR_P/I^2R_P \underset{RP}{\otimes} L$. Thus $\hat{d}$ is injective. Inversely suppose $\hat{d}$ is injective and let $\{f_i \in I\}_{1 \leq i \leq r}$ be such that $\{\bar{f_i} \otimes 1\}$ form a basis of $ IR_P/I^2R_P \underset{RP}{\otimes} L$. Then $\{F_i\}$ form a basis of Im$(\hat{d})$. By fact b), $\widehat{\Omega}_{R/V}$ is a free R-module of rank$(n-1)$ and the $V$-derivations of $R$ generate its dual. The fact that $\{F_i\}$ are linearly independent implies the existence of $V$-derivations $D_i$ of $R$ into itself such that the matrix $(<F_i, D_j \otimes 1>)$ is invertible as stated in ii). 

\smallskip

$\underline {i) \Rightarrow ii)}$. We need to show that i) implies that $\hat{d}$ is injective. 

\smallskip
Note that $\hat{d}$ is the composite of

$I/I^2 \otimes L \overset{\alpha}{\rightarrow} P/(P^2+pR) \otimes L \overset{\beta}{\rightarrow} \hat{\Omega}_{R/V} \otimes L$, where $\beta$ is induced by the standard differential map $"d"$.

Since $R_P$ is a regular local ring and ${I R_P, p}$ form a part of a regular system of parameters of $R_P$, it follows that the natural injection $\alpha$ is injective. Hence it is enough to show that

$P/(P^2+pR) \otimes L \overset{\beta}{\rightarrow} \widehat{\Omega}_{R/V} \otimes L ..............(15)$ 

is injective. 

We have an exact sequence

$P/(P^2+pR) \rightarrow \Omega_{R/V} \otimes R/P \rightarrow \Omega_{(R/P)/V} \rightarrow 0$.

Since $(R, m)$ is a Zarishki ring, it follows from (20.7.20, [5]) that

$$P/(P^2+pR) \rightarrow \widehat{\Omega}_{R/V} \otimes R/P \rightarrow \widehat{\Omega}_{(R/P)/V} \rightarrow 0$$ 

is exact. And hence   
$$P/(P^2+pR) \otimes L \overset{\beta}{\rightarrow} \widehat{\Omega}_{R/V} \otimes L \rightarrow \widehat{\Omega}_{(R/P)/V} \otimes L \rightarrow 0 ......(16)$$
is exact. 

It follows from (theorem 3.6, [3]) that, since $p \in P$, rank$(\widehat{\Omega}_{(R/P)/V} \otimes L) = \dim R/P$ and $\widehat{\Omega}_{R/V} \otimes L$ is of rank$(n-1)$. Moreover, since dim$P/(P^2+pR) \otimes L=$height$P-1$, it follows from (16) that 
$$0 \rightarrow P/(P^2+pR) \otimes L \overset{\beta}{\rightarrow} \widehat{\Omega}_{R/V} \otimes L \rightarrow \widehat{\Omega}_{(R/P)/V} \otimes L \rightarrow 0$$ 
is exact. 

\smallskip
Hence our proof is complete.

\smallskip
\textbf{4.3.} Next we extend Suzuki's theorem to the formal set-up ([10]).

\smallskip
\textbf{Suzuki's theorem.} Let $R$ be a complete intersection locality over a field $k$. Then $R$ is normal if and only if $\Omega_{R/k*}$ is torsion-free, where $k* = k$ if characteristic of $k = 0$ or a differential constant field of $R$ if the characteristic is positive.

Our theorem is the following.

\smallskip
\textbf{Theorem.} Let $C = R/I$, where $(R,m)$ is either unramified or equicharacteristic complete regular local ring with coefficient ring $V$ (field $k$) such that $V/pV(k)$ is perfect and $I$ is a complete intersection prime ideal of $R$. Let $\widehat{\Omega}_{C/V}$, $\widehat{\Omega}_{C/k}$ denote the m-adic completions of $\Omega_{C/V}$, $\Omega_{C/k}$ respectively. We have the following.

\smallskip
i) For mixed-characteristic assume that if $p \neq 0 $ in $C$, then every prime ideal $P$ of height 1 in $C$ containing $p$ is such that $p \notin P^2C_P$. If $C$ is normal then $\widehat \Omega_{C/V}$ and $\widehat \Omega_{C/k}$ are torsion-free over $C$.

\smallskip
ii) If $\widehat{\Omega}_{C/V}$ or $\widehat{\Omega}_{C/k}$ is torsion-free then $C$ is normal.

\smallskip
\textbf{Proof.} If $p \in I$, our assertion reduces to the equicharacteristic case. We here prove the mixed characteristic case; the equicharacteristic case would follow by similar arguments. Our goal is to reduce this formal case to a situation where arguments similar to those given by Suzuki ([10]) would work. Let $R = V[[X_1, ..., X_{n-1}$]] and let $I = (f_1, ..., f_r)$. In order to attain such reduction, for part i) our main tool is theorem 4.2, in particular the observation (14) in the proof; and for part ii) several observations from ([3]) and ([5]).

\smallskip
i) Let $d:R \rightarrow \Omega_{R/V}$ be the usual derivation. Let $g(\neq 0) \in R$; then $g = \lim g_t$ where $g_t$ is the sum of terms of $g$ of degree $\leq t$. We have $dg = \lim dg_t$; $dg_t = {\overset{n-1}{\underset{i=1}{\Sigma}}} (\partial g_t/ \partial X_i) dX_i$, $\lim dg_t = {\overset{n-1}{\underset{i=1}{\Sigma}}} \lim (\partial g_t/ \partial X_i) dX_i = {\overset{n-1}{\underset{i=1}{\Sigma}}} (\partial g/ \partial X_i) dX_i$ (by continuity). Since $\Omega_{R/V}$ is not Hausdorff, $dg = \Sigma (\partial g/ \partial X_i)dX_i + \mu$ where $\mu \in \cap m^t \Omega_{R/V}$. let us recall that the following sequence is split exact (lemma in theorem 3.6, [3]):
$$0 \rightarrow \cap m^t \Omega_{R/V} \rightarrow \Omega_{R/V} \overset{\eta}\rightarrow \widehat{\Omega}_{R/V} (= {\overset{n-1} {\underset{1}{\oplus}}} R dX_i) \rightarrow 0$$
We also denote by $d$ the composite of $R \rightarrow \Omega_{R/V} \rightarrow \widehat{\Omega}_{R/V}$; then, for any $f \in R, df=\Sigma (\partial f/ \partial X_i) dX_i$. Since $(R,m)$ is a Zariski ring, it can be checked easily that the exact sequence 
$$I/I^2 \rightarrow \Omega_{R/V} \otimes C \xrightarrow{^\alpha} \Omega_{C/V} \rightarrow 0$$ 
gives rise to the following exact sequence.
$$I/I^2 \xrightarrow{^{\bar d}} \widehat{\Omega}_{R/V} \otimes C \xrightarrow{^{\widehat{\alpha}}} \widehat{\Omega}_{C/V} \rightarrow 0 ...........(17)$$
\smallskip
First we assume that $C$ is normal. If possible, let $w \in \widehat{\Omega}_{C/V}$ and $c(\neq 0) \in C$ be such that $cw = 0$. Then it follows from (17) that 
$$c {\overset{n-1}{\underset{i=1}{\Sigma}}} \lambda_i \bar{d} X_i = {\overset{r}{\underset{j=1}{\Sigma}}} \nu_j \bar{d} f_j ..........(18)$$
where  $\hat{\alpha} (\Sigma \lambda_i \bar{d} X_i) = w, \bar{d} = d \otimes Id_C,\ \lambda_i, \nu_j$$\in C$; hence
$$c \Sigma \lambda_i \bar{d} X_i = {\overset{d}{\underset{j=1}{\Sigma}}} \nu_j {\overset{n-1}{\underset{i=1}{\Sigma}}} (\partial f_j/ \partial X_i) \bar{d} X_i ...........(19)$$
Since  $\bar{d} X_1, ..., \bar{d} X_{n-1}$ are linearly independent in $\widehat{\Omega}_{R/V}\otimes C$, we have 
$$\lambda_i = 1/c {\overset{d}{\underset{j=1}{\Sigma}}} \nu_j (\partial f_j/ \partial X_i) ...........(20)$$.
We need to show that $\nu_j/c \in C, 1\leq j \leq r$. Then from (18) it would follow that $w=0 \in \widehat{\Omega}_{C/V}$.

\smallskip
Let $P$ be a prime ideal of height $1$ of $C$. Then $C_P$ is a regular local ring. By our assumption $C_P$ is formally smooth over $V$ Let $L=k(P)$-the field of fractions of $C_P$. By the proof of proposition 4.2, $C_P$ is formally smooth over $V$ if and only if 
$$I/I^2 \otimes L \xrightarrow{^{\bar{d_L}}} \widehat{\Omega}_{R/V} \otimes L$$ 
is injective, i.e.
$$I/I^2 \otimes L \xrightarrow{^{\bar{d_L}}} {\overset{n-1}{\underset{i=1}{\oplus}}} L d X_i ..........(21)$$ 
is injective. \\

Due to equation (20), this implies that certain $r$ x $r$ minor of $(\partial f_j/ \partial X_i)$ ${1 \leq j \leq r, {1 \leq i \leq {n-1}}}$ is invertible and hence $\nu_j/c \in C_P$ for $1 \leq j \leq d$. Since $C$ is normal, $C = \cap C_P$ and thus $\nu_j/c \in C$ for $1 \leq j \leq d$.  \\

ii) $R$ is geometrically regular over $V$. Hence $\Gamma_{R/V}=0$ (theorem. 2.3 [3]). We have an exact sequence ([11])
$$0 \rightarrow \Gamma_{C/V} \rightarrow I/I^2 \xrightarrow{\bar{d}} \Omega_{R/V} \otimes C \rightarrow \Omega_{C/V} \rightarrow 0 ....(22)$$
Since $R(C)$ is complete and $I$ is a complete intersection ideal of codimension $r$, $I/I^2$ is a finitely generated free complete $C$-module of rank $r$. Hence, $\Gamma_{C/V}$ and im$\bar{d}$ are also finitely generated complete $C$-modules. taking the m-adic completion we obtain an exact sequence (20.7.20, [5])
$$0 \rightarrow N \rightarrow I/I^2 \xrightarrow{\bar{d}} \widehat{\Omega}_{R/V} \otimes C \rightarrow \widehat{\Omega}_{C/V} \rightarrow 0 .....(23)$$
where $\Gamma_{C/V} \subset N$ (we are not changing the notation for $\bar{d}$). \\

Since the $\rank_R \hat{\Omega}_{R/V}=$dim$R-1$ and $\rank_C \hat{\Omega}_{C/V}=$dim$C-1$ (21.9.5 [5], theorem 3.6, [3]), via rank calculation it follows that $N=0$ and hence $\Gamma_{C/V}=0$. Thus 
$$0 \rightarrow I/I^2 \xrightarrow{\bar{d}} \widehat{\Omega}_{R/V} \otimes C \rightarrow \widehat{\Omega}_{C/V} \rightarrow 0 ....(24)$$
is exact. \\

Let us assume that $C$ is not normal. Let $\bar{C}$ denote the normalization of $C$ in $Q(C)$ and let $J=(C:\bar{C})=$ann${_C} \bar{C}/C$. Let $P$ be a prime ideal of height$1$ in $C$ containing $J$. Then $C_P$ is not regular. By proposition 4.2, rank jac ($f_1, ..., f_r$) mod $P<r$.  The arguments for the rest of the proof is similar to that given by Suzuki for his proof of the affine case. There exist $b_1, ..., b_r$ of $C$, at least one of which, say $b_1$, is not in $P$ such that ${\overset{r}{\underset{i=1}{\Sigma}}} b_i \bar{d}f_i=0$ in $\widehat{\Omega}_{R/V} \otimes C/P$. \\

Let $\tilde{J}=J:P$. There exist an element $j \in \tilde{J} \bar{C}$ such that $j \notin C_P ...............................(25)$ \\

Let $j=c/g$, $c, g \in C$. We have 

$j \Sigma b_i \bar{d} f_i \in j P \widehat{\Omega}_{R/V} \otimes C \subset J \bar{C} \widehat{\Omega}_{R/V} \otimes C \subset \widehat{\Omega}_{R/V} \otimes C$ \\

It follows from above and (24) that $j \Sigma b_i \bar{d} f_i = \Sigma c_i \bar{d} X_i, c_i \in C$

i.e. $c \Sigma b_i \bar{d} f_i = g \Sigma c_i \bar{d} X_i$ in $\widehat{\Omega}_{R/V} \otimes C .....................(26)$ \\

From (24) it follows that im$(g \Sigma c_i \bar{d} X_i) = 0$ in $\widehat{\Omega}_{C/V}$. If im$(\Sigma c_i \bar{d} X_i) \neq 0$ in $\widehat{\Omega}_{C/V}$, then this is a torsion element in $\widehat{\Omega}_{C/V}$. If im$(\Sigma c_i \bar{d} X_i)=0$ in $\widehat{\Omega}_{C/V}$, then it follows from (24) that there are $r$ elements $s_1, ..., s_r \in C$ such that $\Sigma c_i \bar{d} X_i = \Sigma s_i \bar{d} f_i$.  \\

Hence it follows from (26) that  \\
$\Sigma (gs_i - cb_i) \bar{d} f_i =0$. 

Due to (24) $\bar{d} f_1, ..., \bar{d} f_r$ are linearly independent in $\widehat{\Omega}_{R/V} \otimes C$. Hence $gs_i = cb_i, 1 \leq i \leq r$, i.e. $j=c/g = s_i/b_i, 1  \leq i \leq r$. Since $b_1 \notin P$, it follows that $j \in C_P$ - a contradiction due to (25).\\

\textbf{Corollary.} Let $(R,m)$ be an unramified or an equicharacteristic complete regular local ring and let $C=R[[X_1, ..., X_n]]/I$ be a module finite integral extension of $R$ such that $I$ is a complete intersection prime ideal. Suppose that a)$Q(C)/Q(R)$ is separable, b)dim$R>1$ if $R$ is mixed characteristic, dim$R>0$ if $R$ is equicharacteristic and c) in mixed-characteristic for any prime ideal $P$ of height$1$ in $C$ containing $p$, $p \notin P^2C/P.$ If $C$ is normal, then $\Omega_{C/V}$ and $\Omega_{C/k}$ are torsion-free over $C$. \\

\textbf{Proof.} We have a commutative diagram of short exact sequences

\[
\begin{CD}
	@. 0 @. 0 \\
	@. @VVV  @VVV \\
	@.\cap m^t \Omega_{R/V} \otimes C @= \cap m^t \Omega_{R/V} \otimes C @. \\
	@. @VVV  @VVV \\
	0 @>>>\Omega_{R/V} \otimes C @>>> \Omega_{C/V} @>>> \Omega_{C/R} @>>> 0 \\
	%  @VVV  @. @VVV \\
	@.  @VVV @VVV  \parallel @. @.(27)\\
	0 @>>>\widehat{\Omega}_{R/V} \otimes C @>>> \widehat{\Omega}_{C/V} @>>> \Omega_{C/R} @>>>0 \\
	@. @VVV @VVV  @. \\
	@. 0 @. 0 @. 
\end{CD}
\]

\smallskip
$\Omega_{C/R}$, being a finitely generated $C$-module, is complete. Since $\Omega_{R/V}$ is a flat $R$-module and $0 \rightarrow \cap m^t \Omega_{R/V} \rightarrow \Omega_{R/V} \rightarrow \widehat{\Omega}_{R/V} \rightarrow 0$ splits, $\cap m^t \Omega_{R/V} \otimes C$ is a flat $C$-module. Since $C$ is an integral domain the assertion follows immediately. \\

\textbf{4.4. Remarks.} 1. The above diagram shows that for any module-finite extension $C$ of a power series ring $R$ such that $Q(C)$ is separably algebraic over $Q(R)$, the torsion sub-modules of $\Omega_{C/V}$ and $\Omega_{C/k}$ imbed into $\widehat{\Omega}_{C/V}$ and $\widehat{\Omega}_{C/k})$ respectively and the corresponding local cohomology modules of $\Omega_{C/V}$ and $\Omega_{C/k}$ are isomorphic to the same of $\widehat{\Omega}_{C/V}$ and $\widehat{\Omega}_{C/k}$ respectively. \\

2. With notations as above let $B$ denote the absolute integral closure of $C(R)$. Let $\tilde{\Omega}_{B/V}=\underrightarrow{\lim} \widehat{\Omega}_{C/V},R \subset C \subset B, C$ as in Remark 1 above, in mixed characteristic/equicharacteristic $0$. Then the torsion sub-module of $\Omega_{B/V}$ imbeds into $\tilde{\Omega}_{B/V}$ and the corresponding local cohomology modules of $\Omega_{B/V}$ and $\tilde{\Omega}_{B/V}$ are isomorphic. This follows from (27) by taking direct limit and the fact that $\cap m^t \Omega_{R/V}$ is $R$-flat and $\cap m^t \Omega_{R/V} \otimes R/pR =0$ \\

\textbf{4.5.} Our final theorem points out a connection between torsion sub-module of $\widehat{\Omega}_{C/V}$ and derivations $\partial/\partial X_i, 1 \leq i \leq n-1$, of $R$. \\

\textbf{Theorem.} Let $C$ be a complete local domain in mixed characteristic and $R=V[[X_1, ..., X_{n-1}]] \hookrightarrow C$ be such that $C$ is a module finite extension of $R$. Let $m$ denote the maximal ideal of $R$ and let $\widehat{\Omega}_{C/V}$ denote the m-adic completion of $\Omega_{C/V}$. Let $T=$ the torsion sub-module of $\widehat{\Omega}_{C/V}$. Then, for any $\omega \in \widehat{\Omega}_{C/V}, \ \omega \in T \iff \tilde{D}_i (\omega)=0$ for $1 \leq i \leq n-1$, where $\tilde{D}_i:{\Omega}_{C/V}$ or $\widehat{\Omega}_{C/V} \rightarrow Q(C)$ is a $C$-linear map induced by the extension of the derivation $\partial/\partial X_i:R \rightarrow R, \partial/\partial X_i (X_j)=\delta_{ij}$, to $D_i:C \rightarrow Q(C)$. \\

\textbf{Proof.} $\Longrightarrow$: is obvious.

$\Longleftarrow$: let $\widehat{d}$ denote the composite of $C \xrightarrow{d} \Omega_{C/V} \rightarrow \widehat{\Omega}_{C/V} (R \xrightarrow{d} \Omega_{R/V} \rightarrow \hat{\Omega}_{R/V})$. It follows from the middle column of (27) that $(\cap m^t \Omega_{R/V}) \otimes C \simeq \cap m^t \Omega_{C/V})$.

It has been pointed out in the beginning of proof of part i),  theorem 4.3 that for $r \in R, \ dr=\Sigma (\partial r/\partial X_i) dX_i + \mu$ where $\mu \in \cap m^t \Omega_{R/V}$ and $\hat{d}r =\Sigma (\partial r/\partial X_i) dX_i$.

Let $\alpha (\neq 0) \in C$. First we want to show that \\

$f'(\alpha) d \alpha = f'(\alpha) \Sigma D_j(\alpha)dX_j+\nu, \nu \in \cap m^t \Omega_{C/V}$,
where $f(X) \in R[X]$ denotes the minimal polynomial for $\alpha$. \\

Let $f(X)=X^h+r_1X^{h-1}+...+r_h, \ f(\alpha)=0$. \\

Then $f'(\alpha) d \alpha + {\overset{h}{\underset{i=1}{\Sigma}}} d(r_i) \alpha^{h-i}=0$.

For $1 \leq j \leq n-1, \ f'(\alpha) D_j (\alpha) + {\overset{h}{\underset{i=1}{\Sigma}}} D_j(r_i) \alpha^{h-i}=0$.

Hence $f'(\alpha) d \alpha - {\overset{n-1}{\underset{j=1}{\Sigma}}} f'(\alpha) D_j (\alpha) dX_j$

$={\overset{h}{\underset{i=1}{\Sigma}}} [d(r_i)-{\overset{n-1}{\underset{j=1}{\Sigma}}} D_j(r_i)dX_j] \alpha^{h-i}$

$=\Sigma \nu_i \alpha^{h-i}, \nu_i =  [d(r_i)-{\overset{n-1}{\underset{j=1}{\Sigma}}} D_j(r_i)dX_j] \in \cap m^t \Omega_{R/V}$.

Hence $f'(\alpha) d \alpha - {\overset{n-1}{\underset{j=1}{\Sigma}}} f'(\alpha) D_j (\alpha) dX_j \in \cap m^t \Omega_{R/V} \otimes C = \cap m^t \Omega_{C/V} ....(28)$

\smallskip
Now let $\omega \in \Omega_{C/V}$; can write $\omega = {\overset{t}{\underset{i=1}{\Sigma}}} \lambda_i d \alpha_i, \lambda_i, \alpha_i \in C$. Let $f_i(X) \in R[X]$ denote the minimal polynomial of $\alpha_i$ and let $f=\pi f'_i(\alpha_i)$. We have $\tilde{D}_j (\omega) = {\overset{t}{\underset{i=1}{\Sigma}}} \lambda_i D_j \alpha_i$.

It follows from (28) that

$f\omega-\Sigma f\tilde{D}_j (\omega)dX_j \in \cap m^t \Omega_{C/V}     .........(29)$

\smallskip
Hence it follows from (29) that if $\tilde{D}_j (\omega)=0$ for $1 \leq j \leq n-1$, then $f\omega \in \cap m^t \Omega_{C/V} \Longrightarrow f \im \omega =0 \in \widehat{\Omega}_{C/V}$, i.e. $\im \omega \in \widehat{\Omega}_{C/V}$ is a torsion element (middle vertical exact sequence in (27)). 

\smallskip
\textbf{4.6. Corollaries.} Notations are as in Remark 2, 4.4. \\

\textbf{1.} $\omega \in T(\tilde{\Omega}_{B/V}) \iff \tilde{D}_i (\omega)=0, 1 \leq i \leq n-1$, where $T(\tilde{\Omega}_{B/V})$ denotes the torsion sub-module of $\tilde{\Omega}_{B/V}$ and $\tilde{D}_i$ denotes the corresponding extension: $\Omega_{B/V}$ or $\tilde{\Omega}_{B/V} \rightarrow Q(B)$ obtained via the extension of the derivation $\partial/\partial X_i, 1 \leq i \leq n-1,$ as mentioned earlier.

\smallskip
This occurs due to the observation that 
$$0 \rightarrow \cap m^t \widehat{\Omega}_{R/V} \otimes B \rightarrow \Omega_{B/V} \rightarrow \tilde{\Omega}_{B/V} \rightarrow 0,$$
obtained by taking direct limit over $C$ of the middle vertical column in (27), is exact.

\smallskip
\textbf{2.} we have a commutative diagram of short exact sequences:

\[
\begin{CD}
	@. @. 0 @. 0 \\
	@. @. @VVV  @VVV \\
	@. @. T(\tilde{\Omega}_{B/V}) @= T(\tilde{\Omega}_{B/V}) @. \\
	@. @. @VVV  @VVV \\
	0 @>>>\widehat{\Omega}_{R/V} \otimes B @>>> \tilde{\Omega}_{B/V} @>>> \Omega_{B/R} @>>> 0 \\
	%  @VVV  @. @VVV\\
	@.   \parallel @. @VVV @VVV @. ......(30)\\
	0 @>>>\widehat{\Omega}_{R/V} \otimes B @>>> W @>>> N @>>> 0 \\
	@. @. @VVV  @VVV \\
	@. @. 0 @. 0
\end{CD}
\]

where $W$ is a $B$-submodule of ${\overset{n-1}{\underset{i=1}{\oplus}}} Q(B)$ and the vertical map $\tilde{D}:\tilde{\Omega}_{B/V} \rightarrow W$ is given by $\tilde{D}=(\tilde{D}_1, ...., \tilde{D}_{n-1})$; moreover, for $i \geq 0$, $\H_m^i (W)=0$ and $\H_m^i(N) \simeq {\overset{n-1}{\underset{i=1}{\oplus}}} \H_m^{i+1} (B)$.

The proof follows easily from the facts: i) $\tilde{\Omega}_{B/V}=\underrightarrow{\lim} \widehat{\Omega}_{C/V}$ and the bottom row of (30) is exact

ii) $W \xrightarrow{p} W$, multiplication by $p$, is an isomorphism, and iii) $\widehat{\Omega}_{R/V} \otimes B \simeq {\overset{n-1}{\underset{1}{\oplus}}}B$.

\smallskip
Note that $W$ is torsion-free and $W/pW = 0$. Thus $W \xrightarrow{p} W$ is an isomorphism and hence $\H_m^i(W)=0$ for $i\geq 0$.

\smallskip
\textbf{3.} $\Omega_{B/R}$ has an image $N$ such that i) $N$ is cohomologically Cohen-Macaulay i.e. $\H_m^i(N)=0$ for $i < n-1$, $\H_m^{n-1}(N) \neq 0$, ii) $(0:p^t)N$ is finitely generated Cohen-Macaulay free $B/p^tB$-module of rank $n-1$ and of dimension $n-1$, iii) $(0:p^t)N$ is a summand of $(0:p^t)\Omega_{B/R}$ and iv) for any $x \neq 0 \in R, (0:x) N$ can be imbedded into a cohomologically Cohen-Macaulay flat $B/xB$-module. Also $\Omega_{B/V}$ can be mapped onto a $B$-module $W$ such that  $\H_m^i(W) = 0$ for $i \geq 0$. 

\smallskip
\textbf{Proof.} Since $\Omega_{B/V} \rightarrow \tilde{\Omega}_{B/V}$ is onto, the second assertion is obvious from corollary 2, (diagram 30).
The last row of (30) implies that $\H_m^i(N) \simeq \H_m^{i+1}({\overset{n-1}{\underset{1}{\oplus}}}B)$, $(0:p^t)N \simeq {\overset{n-1}{\underset{1}{\oplus}}}B/p^tB$ and the last column of (30) implies that $\H_m^{n-1}(\Omega_{B/R}) \rightarrow \H_m^{n-1}(N) \rightarrow 0$ is exact. Due to: a) theorem 5.1 and corollary 5.10 in ([2]), b) $x$ is a non-zero divisor on $W$ and c) $\Omega_{R/V}$ is $R$-flat, our assertion follows from the last row of (30).

\smallskip
\textbf{4.7. Remarks.} 1) Corollary 2 above provides a new proof of the fact  that $H_m^{n-1}(\Omega_{B/R}) \neq 0$ (Theorem 4.1, [3]). Last row of (30) implies that $H_m^{n-1} (N) \simeq {\overset{n-1}{\underset{1}{\oplus}}} H_m^n(B)$ and last column of (30) implies that $H_m^{n-1}(\Omega_{B/R}) \rightarrow H_m^{n-1} (N)$ is onto. Hence, $H_m^{n-1}(\Omega_{B/R}) \neq 0$ ([1], [3]).

\smallskip
2) Using theorem 5.1 in ([2]) it follows from proposition 4.2 in ([3]) that $\H_m^i (\Omega_{B/V}) \simeq \H_m^i (\Omega_{B/R})$ for $i \leq n-2$, and from corollary 1 (3.4) that a) $\H_m^i(\Omega_{B/A)} \simeq \H_m^i(\Omega_{B/AR})$, for $i \leq n-2$; b)  $\H_m^i(\Omega_{B/V)} \simeq \H_m^i(\Omega_{B/A)}$, for $i \leq n-3$ and injectivity of maps at the next level. The crucial question is whether $\H_m^{n-1}(\Omega_{A/V} \otimes B) \simeq \H_m^n(B) \rightarrow \H_m^{n-1}(\Omega_{B/V})$ is non-null.

\smallskip
3) B. Bhatt kindly shared with me his approach for providing another new proof of the fact that $H_m^{n-1}(\Omega_{B/R}) \neq 0$ (theorem 4.1, [3]). He used Cohen-Macaulayness of the $m$-adic completion $\hat{B}$ along with facts from derived categories and perfectoids to accomplish his proof.

\smallskip
4) In equicharacteristic $0$ same proof works for corresponding statements for theorem 4.5, corollary 1 and corollary 2 in 4.6 without the assertion on local cohomology modules. In positive characteristic the same arguments work for corresponding statement for theorem 4.5 with the assumption that $Q(C)$ is separably algebraic over $Q(R)$.

\smallskip
5) Let $T$ denote the torsion sub-module of $\Omega_{B/V}$, let $U$ = Coker ($T\rightarrow \Omega_{B/V}$) and $L$ = Coker ($T\rightarrow \Omega_{B/R}$). Then, as in cororollary 2 and corollary 3 in 4.6 it follows that $H_m^i(U)=0$ for $i\geq 0$; $H_m^i (L)=0$ for $i\neq n-1$ and $H_m^{n-1} (L) \simeq H_m^{n-1}({\underset{1}{\overset{n-1}{\oplus}}}B)$ and for $t>0, (0:p^t)L \simeq {\underset{1}{\overset{n-1}{\oplus}}}B/p^tB$. The above isomorphisms are obtained from the exact sequence:
$$0\rightarrow \Omega_{R/V} \otimes B \rightarrow U \rightarrow L \rightarrow 0.$$
Thus $L$ shares the same properties with $N$ in corollary 3, 4.6 above. \\
Moreover, with $N$ as in corollary 2 above, there exists an onto map $\eta:L\rightarrow N$ such that $H_m^i(ker \ \eta)=0$ for $i\geq 0$.

\end{document}